\documentclass[11pt,letterpaper]{article}
\usepackage[margin=1in]{geometry}
\usepackage[T1]{fontenc}
\usepackage{amsmath,amssymb,amsthm}
\usepackage{xcolor}
\usepackage{setspace}
\usepackage{microtype}
\usepackage{fancyhdr}
\usepackage{needspace}
\usepackage{tikz}
\usepackage{float}
\definecolor{inkpurple}{RGB}{0,0,0}
\definecolor{inkred}{RGB}{232,54,35}
\definecolor{inkblue}{RGB}{23,105,237}
\usepackage[hidelinks]{hyperref}

\newcommand{\ulinehand}[1]{\textbf{#1}}
\newcommand{\Tq}{\mathcal T_q}
\theoremstyle{plain}
\newtheorem{theorem}{Theorem}
\newtheorem{proposition}{Proposition}
\newtheoremstyle{purplelemma}{3pt}{3pt}{\itshape}{0pt}{\color{inkpurple}\bfseries}{.}{.5em}{}
\theoremstyle{purplelemma}
\newtheorem{lemma}{Lemma}
\newtheorem{corollary}[theorem]{Corollary}
\theoremstyle{definition}
\newtheorem{definition}[theorem]{Definition}
\theoremstyle{remark}
\newtheorem{remark}{Remark}
\renewcommand{\thesection}{\arabic{section}.}

\renewcommand{\theHequation}{\thesection.\arabic{equation}}
\newcommand{\pagebreakhand}{\par\bigskip}
\begin{document}
\normalsize

\vspace*{0pt}
\begin{center}
{\LARGE Counterexamples to the Bilu--Linial signing conjecture in every degree $d\ge3$ and an improved upper bound}\par
\vspace{0.5em}
{\large Jonatan Beiruty}\par
\vspace{0.25em}
\vspace{0.75em}
\vspace{0.75em}\rule{0.96\linewidth}{0.6pt}
\end{center}
\begin{abstract}
For every integer $d\ge3$, we construct a finite $d$-regular graph $G$ such
that every edge signing $s$ satisfies
$\|A_s\|>2\sqrt{d-1}$. This disproves the Bilu--Linial signing conjecture
in every degree. We also prove that every finite simple graph of maximum
degree at most $d$ has an edge signing satisfying
\[
\|A_s\|<(1+\sqrt2)\sqrt{d-1}
\approx2.414213\ldots\sqrt{d-1}.
\]
Our proof uses direct, elementary linear-algebraic techniques and improves
the best known universal coefficient for two-sided graph signing bounds.
Using the same proof technique, we also recover the celebrated bound
$2\sqrt{d-1}$ for bipartite graphs.
\end{abstract}

\tableofcontents

\vspace{0.3cm}\rule{0.43\linewidth}{0.6pt}

\setcounter{section}{-1}\section[Introduction]{Introduction}

Let $G=(V,E)$ be a simple undirected graph of maximum degree at most $d$, where $d\ge3$. We work in the Hilbert space $\mathcal H=\ell^2(V)$ over $\mathbb R$. An edge signing is a function $s:E\to\{-1,+1\}$. Its signed adjacency operator $A_s:\mathcal H\to\mathcal H$ is defined by
\[
(A_sf)(v)=\sum_{w\sim v}s(vw)f(w).
\]
This operator is bounded and is self-adjoint (see \hyperref[lem:signing-facts]{Lemma~1}). When $G$ is finite, the spectral theorem gives $\|A_s\|=\max_{\lambda\in\operatorname{spec}(A_s)}|\lambda|$.

A signing $s$ determines a $2$-lift of $G$, whose eigenvalues consist of the eigenvalues of $G$ together with those of $A_s$. Bilu and Linial stated their conjecture in terms of $2$-lifts~\cite{BiluLinial}:
\begin{center}
\fbox{\begin{minipage}{0.88\linewidth}
\begin{quote}
We conjecture that every $d$-regular graph has a $2$-lift such that all new eigenvalues are in the range $[-2\sqrt{d-1},2\sqrt{d-1}]$.
\end{quote}
\end{minipage}}
\end{center}
Since the new eigenvalues of the $2$-lift are precisely those of $A_s$, this is equivalent to asking whether every finite simple $d$-regular graph has a signing $s$ satisfying
\[
\|A_s\|\le 2\sqrt{d-1}.
\]

\paragraph{Timeline of main results.}
The timeline below records the main developments in signing bounds and counterexamples. The years indicate the first preprint release.
\begin{itemize}
  \item \textbf{2003 --- Bilu and Linial.} Proved a signing bound of order $\sqrt{d\log^3 d}$ for graphs of maximum degree at most $d$, and proposed the conjectured $2\sqrt{d-1}$ bound for regular graphs~\cite{BiluLinial}.
  \item \textbf{2013 --- Marcus, Spielman, and Srivastava.} Proved the conjectured bound for bipartite regular graphs. Their result controls the largest signed adjacency eigenvalue; bipartite spectral symmetry gives the two-sided norm bound~\cite{MSS}.
  \item \textbf{2017 --- Ravichandran and Srivastava.} Obtained the bound $2\sqrt{2(d-1)}$ for regular graphs~\cite{RS2017}.
  \item \textbf{July 2026 --- Huang.} Obtained the degree-sensitive bound $4\sqrt{\lceil d/2\rceil-1}$ for graphs of maximum degree at most $d$~\cite{Huang2026}.
  \item \textbf{September 2026 --- Lin and Zhou.} Obtained the bound $\frac{3+\sqrt5}{2}\sqrt{d-1}\approx2.618\sqrt{d-1}$ for regular graphs~\cite{LinZhou2026}.
  \item \textbf{September 2026 --- Xu.} Constructed a finite simple $3$-regular graph for which every signing has norm greater than $2\sqrt2$, disproving the conjecture in the cubic case~\cite{Xu2026}.
\end{itemize}

In this paper, we prove the following two theorems. The first extends Xu's counterexample result from degree $3$ to every degree. The second improves the known universal upper bounds for the two-sided norm of signed adjacency matrices of graphs with maximum degree at most $d$.

\begin{theorem}[Counterexamples in every degree]\label{thm:intro-counterexamples}
For every integer $d\ge3$, there exists a finite simple $d$-regular graph $G$ such that every edge signing $s:E(G)\to\{-1,+1\}$ satisfies
\[
\|A_s(G)\|>2\sqrt{d-1}.
\]
\end{theorem}

\begin{theorem}[Existence of signings with small spectral norm]\label{thm:intro-small-norm}
For every integer $d\ge3$ and every finite simple graph $G$ of maximum degree at most $d$, there exists an edge signing $s$ such that
\[
\|A_s\|<(1+\sqrt2)\sqrt{d-1}.
\]
\end{theorem}

The proof of the second theorem uses direct, elementary linear-algebraic techniques.

\paragraph{Use of AI tools.} This work was developed through extensive
discussions with ChatGPT (GPT-6.1 Sol), which were central to developing the
main proof ideas and assisting with revisions and corrections. I prepared the
final manuscript and take full responsibility for its mathematical content.

\paragraph{Acknowledgments.} I thank Pesenti Lucas for introducing me to
these problems and for many useful discussions. I also thank Daniel Meriaz
for going over the draft and for his useful feedback and discussions.

\section{The counterexample}

Xu recently constructed a finite cubic counterexample to the Bilu--Linial conjecture~\cite{Xu2026}. We extend this result to every degree $d\ge3$.

Our proof proceeds through an infinite counterexample, which we denote by $H_{d,M}$, and then uses compactness to establish nonconstructively the existence of a finite counterexample as an induced subgraph of the infinite graph. The infinite $d$-regular tree has adjacency norm exactly $2\sqrt{d-1}$, while the Bilu--Linial bound is known to hold for bipartite regular graphs by the result of Marcus, Spielman, and Srivastava~\cite{MSS}. This motivates our construction: we add the simplest odd cycles, triangles, to a tree-like graph and ask whether they can force every signing above the conjectured threshold. 

\pagebreakhand
\vspace*{0.0cm}
\subsection[The construction of the graph $H_{d,M}$]{The construction of the graph $H_{d,M}$}

Fix $d\ge3$ and put $q=d-1$. Let $\Tq$ denote the infinite rooted tree in
which every vertex has exactly $q$ children. Its root has degree $q$
and every other vertex has degree $q+1=d$.

\begin{definition}[The graph $H_{d,M}$]
For a positive integer $M$, define $H_{d,M}$ as follows.
\begin{enumerate}
\item Start with a finite complete rooted $q$-ary tree of height $M$:
the root is at level $0$, and every vertex at a level below $M$
has $q$ children. Label its $q^M$ leaves $o_1,\ldots,o_{q^M}$.
We denote by $U$ the upper tree, consisting of the vertices at levels
$0,\ldots,M-1$ of this rooted tree.\label{def:upper-tree}
\item For each $1\le k\le q^M$, take a disjoint triangle with vertices
$t_k,l_k,r_k$ (top, left, and right), and add the edge $o_kt_k$.
We call the graph obtained after these first two steps the
\emph{finite core}, denoted by $C$.\label{def:finite-core}
\item At every vertex $v$ of this finite core except the root, attach one fresh copy of
$\Tq$, until the degree of the vertex is $d$, joining the vertex to the root of that
copy by a single edge. We call each such attached copy a \emph{branch}
of the vertex $v$, and denote its root by $z_v$. When several branches
are attached to $v$, we distinguish their roots as $z_{v,i}$,
$1\le i\le m(v)$, where $m(v)$ is defined below.
\end{enumerate}

The graph is connected, simple, and countably infinite.
Every vertex except the root has degree $d$; the root has degree $d-1$.
It has precisely $q^M$ cycles, all triangles.
\label{def:branch-count}
For a vertex $v$ of the finite core, let $m(v)$ be the number of copies
of $\Tq$ attached to $v$. For every nonroot vertex, $m(v)=d-\deg_C(v)$.
Explicitly,
\[
 m(v)=
 \begin{cases}
 0,&v\text{ is the finite tree's root},\\
 q-1,&v=o_k,\ l_k,\text{ or }r_k\quad(1\le k\le q^M),\\
 q-2,&v=t_k\quad(1\le k\le q^M),\\
 0,&v\text{ is a finite-tree vertex at a level }1,\ldots,M-1.
 \end{cases}
\]
\end{definition}

\begin{figure}[H]
\centering
\tikzset{
 corevertex/.style={circle,fill=black,inner sep=1.5pt},
 coreedge/.style={line width=0.65pt},
 triangleedge/.style={blue!65!black,line width=0.85pt},
 attachededge/.style={teal!65!black,line width=0.65pt},
 treecircle/.style={circle,draw=teal!65!black,densely dotted,
   line width=0.9pt,minimum size=0.66cm,inner sep=1pt,
   text=teal!65!black,font=\small}
}
\begin{minipage}[t]{0.30\textwidth}
\vspace{0pt}\centering
\begin{tikzpicture}[x=0.82cm,y=0.95cm,every node/.style={font=\small}]
 \coordinate (root) at (0,2.65);
 \coordinate (left) at (-1.1,1.6);
 \coordinate (right) at (1.1,1.6);
 \draw[coreedge] (root)--(left) (root)--(right);
 \draw[coreedge] (left)--(-1.65,0.65) (left)--(-0.55,0.65);
 \draw[coreedge] (right)--(0.55,0.65) (right)--(1.65,0.65);
 \foreach \p in {root,left,right} \node[corevertex] at (\p) {};
 \node[above] at (root) {$z$};
 \foreach \x in {-1.65,-0.55,0.55,1.65} {
   \node[corevertex] at (\x,0.65) {};
   \draw[coreedge] (\x,0.65)--({\x-0.19},0.12)
                  (\x,0.65)--({\x+0.19},0.12);
   \node at ({\x-0.19},-0.15) {$\vdots$};
   \node at ({\x+0.19},-0.15) {$\vdots$};
 }
\end{tikzpicture}
\par\smallskip
{\small (a) $\Tq$, shown for $q=2$.}
\end{minipage}\hfill
\begin{minipage}[t]{0.68\textwidth}
\vspace{0pt}\centering
\begin{tikzpicture}[x=1.00cm,y=0.88cm,every node/.style={font=\small}]
 \coordinate (root) at (0,3.2);
 \coordinate (left) at (-2.2,2.1);
 \coordinate (right) at (2.2,2.1);
 \draw[coreedge] (root)--(left) (root)--(right);
 \draw[coreedge] (left)--(-3.3,1.0) (left)--(-1.1,1.0);
 \draw[coreedge] (right)--(1.1,1.0) (right)--(3.3,1.0);
 \foreach \p in {root,left,right} \node[corevertex] at (\p) {};
 \node[above] at (root) {root};
 \foreach \k/\x in {1/-3.3,2/-1.1,3/1.1,4/3.3} {
   \coordinate (o\k) at (\x,1.0);
   \coordinate (t\k) at (\x,0.2);
   \coordinate (l\k) at ({\x-0.38},-0.5);
   \coordinate (r\k) at ({\x+0.38},-0.5);
   \draw[coreedge] (o\k)--(t\k);
   \draw[triangleedge] (t\k)--(l\k)--(r\k)--cycle;
   \foreach \letter in {o,t,l,r} \node[corevertex] at (\letter\k) {};
   \ifodd\k\relax
     \node[right,xshift=2pt,yshift=-2pt] at (o\k) {$o_{\k}$};
   \else
     \node[right] at (o\k) {$o_{\k}$};
   \fi
   \node[right] at (t\k) {$t_{\k}$};
   \node[below left] at (l\k) {$l_{\k}$};
   \node[below right] at (r\k) {$r_{\k}$};
 }
\end{tikzpicture}
\par\smallskip
{\small (b) The finite core $C$ for $q=2$, $M=2$.}
\end{minipage}

\par\medskip
\begin{tikzpicture}[x=1.20cm,y=0.96cm,every node/.style={font=\small}]
 \coordinate (root) at (0,3.3);
 \coordinate (left) at (-2.2,2.3);
 \coordinate (right) at (2.2,2.3);
 \draw[coreedge] (root)--(left) (root)--(right);
 \draw[coreedge] (left)--(-3.3,1.2) (left)--(-1.1,1.2);
 \draw[coreedge] (right)--(1.1,1.2) (right)--(3.3,1.2);
 \foreach \p in {root,left,right} \node[corevertex] at (\p) {};
 \node[right,xshift={5pt+0.2cm}] at (root) {root};
 \foreach \k/\x in {1/-3.3,2/-1.1,3/1.1,4/3.3} {
   \coordinate (o\k) at (\x,1.2);
   \coordinate (t\k) at (\x,0.25);
   \coordinate (l\k) at ({\x-0.38},-0.55);
   \coordinate (r\k) at ({\x+0.38},-0.55);
   \draw[coreedge] (o\k)--(t\k);
   \draw[triangleedge] (t\k)--(l\k)--(r\k)--cycle;
   \foreach \letter in {o,t,l,r} \node[corevertex] at (\letter\k) {};
   \ifodd\k\relax
     \node[right,xshift=2pt,yshift=-2pt] at (o\k) {$o_{\k}$};
   \else
     \node[above right] at (o\k) {$o_{\k}$};
   \fi
   \node[right] at (t\k) {$t_{\k}$};
   \node[above left] at (l\k) {$l_{\k}$};
   \node[above right] at (r\k) {$r_{\k}$};
   \node[treecircle] (ob\k) at ({\x-0.78},1.2) {$\Tq$};
   \node[treecircle] (lb\k) at ({\x-0.75},-1.45) {$\Tq$};
   \node[treecircle] (rb\k) at ({\x+0.75},-1.45) {$\Tq$};
   \draw[attachededge] (o\k)--(ob\k) (l\k)--(lb\k) (r\k)--(rb\k);
 }
\end{tikzpicture}
\par\smallskip
{\small (c) The same finite core with its attachments: $H_{3,2}$.}

\caption{The graph $H_{3,2}$ constructed in 3 steps.}
\label{fig:construction}
\end{figure}
\subsection{proof $H_{d,M}$ is a counterexample for large $M$.}

To prove that $H_{d,M}$ is a counterexample for large $M$, we first establish a simple lemma about signed adjacency operators. The lemma gives the norm bounds we need and lets us simplify arbitrary edge signings.

\vspace{0.35cm}
\begin{lemma}[General signing facts]\label{lem:signing-facts}
Let $G$ be a finite or countably infinite simple graph of maximum
degree $d<\infty$, let $s:E(G)\to\{-1,1\}$, and write
$A=A_s(G)$ on $\ell^2(V(G))$. A cycle's \emph{flux} is
$\sigma_s(C)=\prod_{e\in E(C)}s(e)$.\label{def:flux}
\begin{enumerate}
\item\textbf{Operator bounds.}
$A$ is bounded and self-adjoint, with $\|A\|\le d$.
For every $W\subseteq V(G)$, let $F=G[W]$ be the induced subgraph on $W$.
Zero extension $\iota:\ell^2(W)\to\ell^2(V(G))$ gives
\[
 A_{s|_{E(F)}}(F)=\iota^*A\iota,\qquad
 \|A_{s|_{E(F)}}(F)\|\le\|A\|.
\]

\item\textbf{Switching and flux.}
For any vertex-signing function $\eta:V(G)\to\{-1,1\}$, define the
switched edge signing by
$s^\eta(uv)=\eta(u)s(uv)\eta(v)$. This switching preserves the
operator norm, $\|A_{s^\eta}\|=\|A_s\|$, and the flux of every cycle.
Thus, in a connected graph, we can choose $\eta$ so that every edge of
a given spanning tree has sign $+1$.
\end{enumerate}
\end{lemma}

\begin{proof}[Proof of Lemma 1]
\textbf{1. Proof of operator bounds.}

There are three claims to prove here: boundedness, self-adjointness,
and the norm bound for every induced subgraph.

For finitely supported $f\in\ell^2(V(G))$, Cauchy--Schwarz gives
\[
\|Af\|^2=\sum_{v\in V(G)}\left|\sum_{w\sim v}s(vw)f(w)\right|^2\le d\sum_{v\in V(G)}\sum_{w\sim v}|f(w)|^2=d\sum_{w\in V(G)}\deg(w)|f(w)|^2\le d^2\|f\|^2.
\]
Since finitely supported functions are dense in $\ell^2(V(G))$, $A$
extends to a bounded operator with $\|A\|\le d$.

The signed adjacency matrix is symmetric, so for finitely supported
$f,g$, $\langle Af,g\rangle=\langle f,Ag\rangle$.
By density this holds for all $f,g\in\ell^2(V(G))$, hence $A$ is
self-adjoint.

Let $W$ be a subset of vertices of $G$. For $F=G[W]$, let
$\iota:\ell^2(W)\to\ell^2(V(G))$ be its zero
extension. It is an isometry, and its adjoint restricts functions to
$W$. Therefore
\[
A_{s|_{E(F)}}(F)=\iota^*A\iota,
\qquad
\|A_{s|_{E(F)}}(F)\|\le\|A\|.
\]

\textbf{2. Proof of switching and flux.} Let $\eta$ be a vertex-signing function $\eta:V(G)\to\{-1,1\}$, and define the linear operator $D_\eta:\ell^2(V(G))\to\ell^2(V(G))$ by $(D_\eta f)(v)=\eta(v)f(v)$.
For every $f\in\ell^2(V(G))$,
\[
\|D_\eta f\|_2^2
=\sum_{v\in V(G)}|\eta(v)f(v)|^2
=\sum_{v\in V(G)}|f(v)|^2
=\|f\|_2^2.
\]
Thus $D_\eta$ is an isometry, and $D_\eta^2=I$. For every $f\in\ell^2(V(G))$ and $v\in V(G)$, we have
\[
(A_{s^\eta}f)(v)
=\sum_{w\sim v}s^\eta(vw)f(w)
=\sum_{w\sim v}\eta(v)s(vw)\eta(w)f(w)
=\eta(v)\bigl(A_s(D_\eta f)\bigr)(v).
\]
Therefore $A_{s^\eta}=D_\eta A_sD_\eta$. Since $D_\eta$ is an isometry, it follows that $\|A_{s^\eta}\|=\|D_\eta A_sD_\eta\|=\|A_s\|$.

For a cycle $C=v_1v_2\cdots v_kv_1$, with $v_{k+1}=v_1$, its flux is
\[
\sigma_{s^\eta}(C)
=\prod_{i=1}^k s^\eta(v_iv_{i+1})
=\prod_{i=1}^k\eta(v_i)s(v_iv_{i+1})\eta(v_{i+1})
=\prod_{i=1}^k s(v_iv_{i+1})
=\sigma_s(C),
\]
because each vertex sign appears twice and $\eta(v_i)^2=1$. Thus
switching preserves the flux of every cycle.

Now let $T$ be a spanning tree of the connected graph $G$, rooted at
$v_0$. Set $\eta(v_0)=1$. Recursively, for each tree edge from a
parent $v$ to a child $w$, set $\eta(w)=\eta(v)s(vw)$.
Then $s^\eta(vw)=\eta(v)s(vw)\eta(w)=1$, so every edge of $T$ has
sign $+1$ after switching.
\end{proof}

\vspace{0.45cm}
We are ready to prove that $H_{d,M}$ is an infinite counterexample.

\vspace{0.25cm}
\begin{proposition}\label{prop:HdM-counterexample}
There exists a sufficiently large $M$ such that $H_{d,M}$ is a counterexample.
\end{proposition}
Before proving this proposition, we will define a useful quantity and lay out the proof strategy.

\begin{definition}[Energy]
Let $G$ be a graph of bounded degree, let $s$ be an edge signing of $G$, and fix $R\in\mathbb R$. For $f\in\ell^2(V(G))$, the \emph{energy} of $f$ is
\[
Q_s(f):=\langle f,(RI-A_s)f\rangle.
\]
\end{definition}

If $\|A_s\|\le R$, then for every $f\in\ell^2(V(G))$,
\[
Q_s(f)=R\|f\|^2-\langle f,A_sf\rangle
\ge (R-\|A_s\|)\|f\|^2\ge0.
\]
Since $A_{-s}=-A_s$, the same bound holds for $Q_{-s}(f)$.

\paragraph{Proof strategy.} We argue by contradiction.
\begin{enumerate}
\item Suppose that, for some sufficiently large $M$ to be specified below, $H_{d,M}$ has a signing $s$ with $\|A_s\|\le R$, where $R=2\sqrt{d-1}$.
\item By the switching and flux claims of Lemma~1, switch $s$ and $-s$ so that both signings are $+1$ on the chosen spanning tree. Their cycle fluxes are preserved. Denote the resulting signings by $s_+$ and $s_-$.
\item Since switching preserves the norm, $\|A_{s_\pm}\|=\|A_s\|\le R$. Therefore, the energies $Q_\pm(f):=Q_{s_\pm}(f)$ are nonnegative for every $f$:
\begingroup
\renewcommand{\theHequation}{\thesection.star}
\begin{equation}\label{eq:nonnegative-switched-energies}
Q_\pm(f):=\langle f,(RI-A_{s_\pm})f\rangle\ge0.
\tag{$*$}
\end{equation}
\endgroup
\Needspace{9\baselineskip}
\item To reach a contradiction, we would like to construct a vector with negative energy. We will in fact construct two ``twin vectors'' $f_+$ and $f_-$ with negative average energy, so at least one of them has negative energy. To find such a pair, we construct a sequence of twin vectors $f_+^\gamma,f_-^\gamma$ whose average energy has a negative limit:
\[
\lim_{\gamma\uparrow1}\frac{Q_+(f_+^\gamma)+Q_-(f_-^\gamma)}{2}<0.
\]
For $\gamma<1$ sufficiently close to $1$, the average is then negative, so at least one energy is negative, contradicting~\eqref{eq:nonnegative-switched-energies}.
\end{enumerate}

For the rest of the proof, put $b=q^{-1/2}=(d-1)^{-1/2}$ and let
\[
T:=H_{d,M}\setminus\{l_kr_k:1\le k\le q^M\}.
\]
By the construction of $H_{d,M}$, $T$ is a spanning tree. For each $k$, let
$\sigma_k$ be the flux of the $k$th triangle in the original signing $s$.
Switch $s$ and $-s$ separately so that every edge of $T$ has sign $+1$,
and denote the resulting signings by $s_+$ and $s_-$. The flux of the $k$th
triangle is preserved by switching; since its two edges other than $l_kr_k$
belong to $T$, this gives
\[
s_+(l_kr_k)=\sigma_k,\qquad s_-(l_kr_k)=-\sigma_k,
\]
while all other edges have sign $+1$. For $\theta\in\{+1,-1\}$, write
$s_\theta=s_+$ if $\theta=+1$ and $s_\theta=s_-$ if $\theta=-1$, and set
$Q_\theta(f):=Q_{s_\theta}(f)$.

\Needspace{5\baselineskip}
\paragraph{Energy decomposition.} Recall that $H_{d,M}$ consists of the finite core $C$ and the attached branches $\Tq^\beta$. For each branch, let $v_\beta\in C$ be its attachment vertex and $z_\beta$ its root. Then, for every signing $s$ and every $f\in\mathcal H$,
\begin{equation}\label{eq:energy-decomposition}
\begin{aligned}
Q_s(f)={}&\underbrace{R\sum_{v\in C}f(v)^2-2\sum_{\{u,v\}\in E(C)}s(uv)f(u)f(v)}_{\text{finite-core contribution}}\\
&+\sum_\beta\underbrace{\left[R\sum_{w\in V(\Tq^\beta)}f(w)^2-2\sum_{\{u,w\}\in E(\Tq^\beta)}s(uw)f(u)f(w)-2s(v_\beta z_\beta)f(v_\beta)f(z_\beta)\right]}_{\text{branch contribution, including its connecting edge}}.
\end{aligned}
\end{equation}
The final term in each branch contribution accounts for the edge joining the branch to the finite core.

We are now ready to prove Proposition~\ref{prop:HdM-counterexample}.

\begin{proof}[Proof of Proposition~\ref{prop:HdM-counterexample}]

Set $R=2\sqrt{d-1}$. Fix $M\ge1$ and suppose that $H_{d,M}$ has a signing $s$ with $\|A_s\|\le R$. Apply the switching claim of Lemma~1 to $s$ and $-s$, and denote the resulting signings by $s_+$ and $s_-$. By~\eqref{eq:nonnegative-switched-energies}, the energies $Q_+(f)$ and $Q_-(f)$ are nonnegative for every $f\in\mathcal H$. We derive a contradiction by constructing a sequence of ``twin vectors'' $f_+^\gamma,f_-^\gamma$ such that
\[
\lim_{\gamma\uparrow1}\frac{Q_+(f_+^\gamma)+Q_-(f_-^\gamma)}2<0.
\]
The calculation below shows that this limit is negative for all sufficiently large $M$. Then, for $\gamma<1$ sufficiently close to $1$, the average energy is negative, so one of the two energies is negative.

\vspace{0.45cm}
We are now ready to define $f_\pm^\gamma$. To keep the notation light, write $f_\pm=f_\pm^\gamma$. We will define $f_\pm$ in 3 steps. The first will be to define $f_\pm$ on the infinite branches as a function of the value of the vertex of each finite branch. Then in the second step we will define the values on each triangle and its attached old leaf as functions of the value at the old leaf's parent. And lastly we will define the remaining values of the finite core on the upper tree $U$, at levels $0,\ldots,M-1$.

\vspace*{0.25cm}
\noindent{\Large\ulinehand{Step 1: values on infinite branches.}} For $0<\gamma<1$, once the values on the finite core are chosen, set
\[
 f_\theta(w)=f_\theta(v)(b\gamma)^j
 \qquad\text{if }w\text{ is at distance }j\ge1\text{ from }v,
\]
for each attached branch and $\theta\in\{+1,-1\}$.

\vspace{0.35cm}
\noindent{\Large\ulinehand{Step 2: values on triangles and old leaves.}} For each old leaf $o_k$, let $p(o_k)$ denote its parent in the finite tree, and put
$a=f_\theta(p(o_k))$, the value at that parent (which will be the same for $f_+$ and $f_-$ and for all $k$). Set
\[
\begin{aligned}
 u_k:=f_\theta(o_k)&=a(b+\theta\sigma_kb^6),\\
 b_k:=f_\theta(t_k)&=a(b^2+\theta\sigma_kb^5),\\
 w_k:=f_\theta(l_k)=f_\theta(r_k)&=a(b^3+\theta\sigma_kb^4).
\end{aligned}
\]

\vspace{0.35cm}
\noindent{\Large\ulinehand{Step 3: values on the remaining of the finite core.}} For each vertex $v\in U$ at level $j=0,\ldots,M-1$, set
\[
 f_+(v)=f_-(v)=(j+1)b^j.
\]

\vspace{0.35cm}
Using the energy decomposition~\eqref{eq:energy-decomposition}, we now calculate the core and branch contributions to $Q_\theta(f_\theta)$.

\noindent\textbf{Infinite-branch contribution.}
For a branch attached at $v_\beta$, direct summation gives
\[
\begin{aligned}
&R\sum_{w\in V(\Tq^\beta)}f_\theta(w)^2
 -2\sum_{\{u,w\}\in E(\Tq^\beta)}f_\theta(u)f_\theta(w)
 -2f_\theta(v_\beta)f_\theta(z_\beta)\\
&=f_\theta(v_\beta)^2\left[
 R\sum_{j=1}^\infty q^{j-1}(b\gamma)^{2j}
 -2\sum_{j=1}^\infty q^j(b\gamma)^{2j+1}
 -2b\gamma\right]\\
&=\left(-b+\eta_\gamma\right)f_\theta(v_\beta)^2,
\qquad \eta_\gamma=b\frac{1-\gamma}{1+\gamma}\longrightarrow0
\quad(\gamma\uparrow1).
\end{aligned}
\]

\medskip\noindent\textbf{Finite-core contribution.}
We first write the finite-core contribution using the values from Step~3:
\[
\begin{aligned}
R\sum_{v\in C}f_\theta(v)^2
&-2\sum_{\{u,v\}\in E(C)}s_\theta(uv)f_\theta(u)f_\theta(v)\\
={}&2b^{-1}\left[
\sum_{j=0}^{M-1}(j+1)^2-
\sum_{j=0}^{M-2}(j+1)(j+2)\right]\\
&+\sum_{k=1}^{q^M}\Bigl[
2b^{-1}(u_k^2+b_k^2+2w_k^2)
-2au_k-2u_kb_k-4b_kw_k-2\theta\sigma_k w_k^2\Bigr].
\end{aligned}
\]
Substituting the Step~2 values and simplifying the sums, we get
\[
R\sum_{v\in C}f_\theta(v)^2
-2\sum_{\{u,v\}\in E(C)}s_\theta(uv)f_\theta(u)f_\theta(v)
=b^{-1}M(M+1)-2M^2b^7
+\frac{M^2}{q^M}\left(6b^4-8b^6\right)\sum_{k=1}^{q^M}\theta\sigma_k.
\]

Adding the finite-core and infinite-branch calculations, we obtain
\begin{equation}\label{eq:combined-core-branch-energy}
\small
Q_\theta(f_\theta)=b^{-1}M(M+1)-2M^2b^7
+\frac{M^2}{q^M}\left(6b^4-8b^6\right)\sum_{k=1}^{q^M}\theta\sigma_k
+\left(-b+\eta_\gamma\right)\sum_\beta f_\theta(v_\beta)^2.
\end{equation}
Taking the average of the two endpoint energies in~\eqref{eq:combined-core-branch-energy}, the terms containing $\theta\sigma_k$ cancel, and we obtain
\[
\frac{Q_+(f_+)+Q_-(f_-)}2
=b^{-1}M(M+1)-2M^2b^7
+\frac12\left(-b+\eta_\gamma\right)
\sum_{\theta=\pm1}\sum_\beta f_\theta(v_\beta)^2.
\]
Recall that $v_\beta$ is the vertex in the finite core to which the infinite branch $\Tq^\beta$ is attached. The vertices in $C$ connected to infinite branches are exactly those with $m(v)>0$. Therefore,
\[
\sum_{\theta=\pm1}\sum_\beta f_\theta(v_\beta)^2
=\sum_{\theta=\pm1}\sum_{v\in C}m(v)f_\theta(v)^2
=2q^M\bigl(f_\theta(p(o_k))\bigr)^2\left(1-b^8-b^{10}-b^{12}\right).
\]
Let $\delta:=b^8-b^{10}-b^{12}=b^8\left(1-b^2-b^4\right)>0$.
Combining this sum with the finite-core contribution gives
\[
\frac{Q_+(f_+)+Q_-(f_-)}2
=b^{-1}M(M+1)
+q^M\bigl(f_\theta(p(o_k))\bigr)^2\left[-b(1+\delta)
+\eta_\gamma\left(1-b^8-b^{10}-b^{12}\right)\right].
\]
Moreover, substituting the common value $f_\theta(p(o_k))=Mb^{M-1}$, we get
\[
\frac{Q_+(f_+)+Q_-(f_-)}2
 =b^{-1}\left[M(M+1)-(1+\delta)M^2\right]
 +b^{-2}\eta_\gamma M^2\left(1-b^8-b^{10}-b^{12}\right).
\]
Since $\eta_\gamma\to0$ as $\gamma\uparrow1$, the formula above gives
\Needspace{7\baselineskip}
\[
\lim_{\gamma\uparrow1}\frac{Q_+(f_+^\gamma)+Q_-(f_-^\gamma)}2
=b^{-1}\left[M(M+1)-(1+\delta)M^2\right]
=b^{-1}(M-\delta M^2).
\]
This is negative for all sufficiently large $M$. For such an $M$ take $\gamma<1$ sufficiently close to 1 so that the average is negative and we are done.

\end{proof}

\begin{remark}
Viewing $H_{d,M}$ as a perturbation of the infinite $d$-regular tree $\mathbb T_d$, we choose the twin vectors $f_+^\gamma$ and $f_-^\gamma$ to resemble the radial test vectors used to show that $\|A_{\mathbb T_d}\|\ge 2\sqrt{d-1}$. Both use the decay $b^j=(d-1)^{-j/2}$. Our vectors are modified near the triangles to account for their cycle fluxes, and along the infinite branches by a factor $\gamma<1$ to ensure square-summability.
\end{remark}

\pagebreakhand
\vspace*{0pt}
\subsection{from an infinite counterexample to a finite one}

We now show that an infinite counterexample implies the existence of a finite one. This works for any infinite counterexample; we do not need to use the special geometry of our infinite counterexample.\par\smallskip
The first lemma shows that any connected infinite counterexample of maximum degree at most $d$ contains a finite connected induced subgraph that is still a counterexample. The second lemma shows that every finite connected graph of maximum degree at most $d$ is an induced subgraph of a finite connected $d$-regular graph. Thus, for this existence question, it is equivalent to consider finite graphs of maximum degree at most $d$ or finite $d$-regular graphs.
\medskip
\begin{lemma}[Compactness]\label{lem:compactness}
Let $H$ be a connected graph of maximum degree at most $d$, and let
$R\ge0$. If $\|A_s(H)\|>R$ for every signing
$s:E(H)\to\{-1,+1\}$, then $H$ has a finite, connected, induced
subgraph $F$ such that $\|A_\tau(F)\|>R$ for every signing
$\tau:E(F)\to\{-1,+1\}$.
\end{lemma}

\vspace{0.45cm}
The key idea of the proof is that each signing has a finitely supported witness (that is, a finitely supported vector $f$ with $|\langle A_s(H)f,f\rangle|>R$). Compactness of the signing space $\{-1,+1\}^{E(H)}$ then lets us find a single finite induced subgraph on which every signing still has norm greater than $R$.

\begin{proof}
Let $\Omega=\{-1,+1\}^{E(H)}$ be the space of edge signings of $H$,
equipped with the product topology, making it a compact topological space.

By assumption, $\|A_s(H)\|>R$ for every $s\in\Omega$. By the operator
bounds in \hyperref[lem:signing-facts]{Lemma~1}, $A_s(H)$ is a bounded
self-adjoint operator. The spectral theorem therefore gives
\[
\|A_s(H)\|=\sup_{\|f\|_2=1}
\bigl|\langle A_s(H)f,f\rangle\bigr|.
\]
Thus, for every $s\in\Omega$, there is a unit vector with quadratic form
of absolute value greater than $R$. By density of the finitely supported
vectors in $\ell^2(V(H))$ and continuity of the quadratic form, we may
choose such a vector to be finitely supported. Denote it by $f_s$; then
\[
\|f_s\|_2=1,
\qquad
\bigl|\langle A_s(H)f_s,f_s\rangle\bigr|>R.
\]
Let $E_s=\{\{u,v\}\in E(H):u,v\in\operatorname{supp}f_s\}$, the
finite set of edges whose endpoints both lie in the support of $f_s$,
and define
\[
U_s=\{s'\in\Omega:s'(e)=s(e)\text{ for every }e\in E_s\}.
\]
Since $U_s$ fixes only finitely many coordinates, it is open in the
product topology, and $s\in U_s$.
Since this holds for every $s\in\Omega$, the sets $U_s$ cover $\Omega$:
$\bigcup_{s\in\Omega}U_s=\Omega$. By compactness of $\Omega$, there is
a finite subcover $U_{s_1},\ldots,U_{s_k}$, so
$\bigcup_{i=1}^k U_{s_i}=\Omega$.

Let $W$ contain the union of the supports of
$f_{s_1},\ldots,f_{s_k}$ and the vertices along finitely many paths in
$H$ connecting all vertices in these supports. Since $H$ is connected,
such paths exist; since there are finitely many support vertices, $W$ is
finite. Define $F=H[W]$, the induced subgraph of $H$ on $W$. Then $F$ is
finite and connected.

We claim that every signing of $F$ has norm greater than $R$. Let
$\tau:E(F)\to\{-1,+1\}$ be any signing, and extend it to a signing
$\widetilde\tau$ of $H$, for example by assigning $+1$ to every edge
outside $E(F)$. Since $U_{s_1},\ldots,U_{s_k}$ cover $\Omega$, there
is some $i$ such that $\widetilde\tau\in U_{s_i}$. Thus
$\widetilde\tau(e)=s_i(e)$ for every $e\in E_{s_i}$. Because
$\operatorname{supp}f_{s_i}\subseteq W$ and $F=H[W]$ is induced,
\[
\bigl|\langle A_\tau(F)f_{s_i}|_W,f_{s_i}|_W\rangle\bigr|
=\bigl|\langle A_{\widetilde\tau}(H)f_{s_i},f_{s_i}\rangle\bigr|
=\bigl|\langle A_{s_i}(H)f_{s_i},f_{s_i}\rangle\bigr|>R.
\]
Since $f_{s_i}|_W$ is a unit vector, this implies
$\|A_\tau(F)\|>R$. As $\tau$ was arbitrary, every signing of $F$
has norm greater than $R$.
\end{proof}

This lemma lets us reduce the $d$-regular case to graphs of maximum
degree at most $d$: it suffices to find a finite simple graph of maximum
degree at most $d$ for which every signing violates the conjectured
bound, since it can be completed to a $d$-regular graph that also
violates the bound.

\vspace{0.8cm}
\begin{lemma}[Regular completion]\label{lem:regular-completion}
Let $G$ be a finite simple graph of maximum degree at most $d$. Then
$G$ is an induced subgraph of a finite simple $d$-regular graph $G'$.
\end{lemma}

\begin{proof}[Proof of Lemma~3]
Choose an even integer $n>d$. For each $v\in V(G)$, let
$r_v=d-\deg_G(v)$, so $0\le r_v\le d<n$. On the cyclic group
$\mathbb Z_n$, choose a simple $r_v$-regular graph $R_v$. Such a graph
can be obtained by joining each $i\in\mathbb Z_n$ to
$i\pm1,\ldots,i\pm r_v/2$ when $r_v$ is even; when $r_v$ is odd,
join it to $i\pm1,\ldots,i\pm(r_v-1)/2$ and to $i+n/2$.

Take $V(G')=V(G)\times\mathbb Z_n$. Include an edge between $(u,i)$
and $(v,i)$ whenever $uv\in E(G)$, and, for each $v\in V(G)$, include
an edge between $(v,i)$ and $(v,j)$ whenever $ij\in E(R_v)$. The
resulting graph is finite and simple. Each vertex $(v,i)$ has
$\deg_G(v)+r_v=d$ neighbors, so $G'$ is $d$-regular.

For any fixed $i\in\mathbb Z_n$, the vertices
$V(G)\times\{i\}$ induce a copy of $G$: the added edges between copies
of a vertex never join two vertices within this layer. Hence $G$ is an
induced subgraph of $G'$.
\end{proof}

Now we are ready to prove Theorem~1, by combining Proposition~2 and
the two lemmas above.

\begin{proof}[Proof of Theorem~1]
Fix $d\ge3$. By Proposition~\ref{prop:HdM-counterexample}, choose
$M$ sufficiently large that $H_{d,M}$ is an infinite counterexample:
every signing $s$ of $H_{d,M}$ satisfies
$\|A_s(H_{d,M})\|>2\sqrt{d-1}$. The graph $H_{d,M}$ has maximum degree
at most $d$. Applying Lemma~\ref{lem:compactness} with
$R=2\sqrt{d-1}$ gives a finite connected induced subgraph $F$ of
$H_{d,M}$ such that every signing $\tau$ of $F$ satisfies
$\|A_\tau(F)\|>2\sqrt{d-1}$. Since $F$ is an induced subgraph of
$H_{d,M}$, its maximum degree is also at most $d$.

By Lemma~\ref{lem:regular-completion}, $F$ is an induced subgraph of a
finite simple $d$-regular graph $G$. Let $s$ be any signing of $G$,
and let $\tau$ be its restriction to $F$. Because $F$ is induced, the
signed adjacency matrix $A_\tau(F)$ is a principal submatrix of
$A_s(G)$. The norm bound for induced subgraphs in
\hyperref[lem:signing-facts]{Lemma~1} therefore
gives
\[
\|A_s(G)\|\ge \|A_\tau(F)\|>2\sqrt{d-1}.
\]
Since $s$ was arbitrary, every signing of $G$ violates the
Bilu--Linial bound. Thus $G$ is the required finite simple
$d$-regular counterexample.
\end{proof}

\begin{remark}
A quantitative refinement of the construction shows that there exists an
absolute constant $c>0$ such that, for every sufficiently large integer
$d$, the counterexample $G_d$ can be chosen so that every edge signing
$s$ satisfies
\[
\|A_s(G_d)\|
\ge 2\sqrt{d-1}+\frac{c}{(d-1)^{15/2}}.
\]
\end{remark}

\begin{remark}
To avoid first constructing an infinite counterexample and then using
compactness to obtain a finite one nonconstructively, one may work
directly with $H'_{d,M}$, obtained from $H_{d,M}$ by replacing each
infinite branch $\mathcal T_q$ with its first $M$ levels. The same
proof, with slightly messier constants, shows that for each fixed
$d\ge3$ and all sufficiently large $M$, every signing $s$ satisfies
\[
\|A_s(H'_{d,M})\|>2\sqrt{d-1}.
\]
The regular-completion lemma then gives a finite $d$-regular
counterexample.
\end{remark}

\pagebreakhand
\section{Existence of Signings with Small Spectral Norm}

Motivated by the determinant weight in the recent algorithm of Jadbabaie, Saberi, and Sra~\cite{JSS2026}, we define a score using that weight on partial signings and show that a maximal-score partial signing can be extended. This yields an existence proof for a signing with small spectral norm.

\setcounter{theorem}{1}
\begin{theorem}\label{thm:small-norm-signing}
Let $G=(V,E)$ be a finite simple undirected graph of maximum degree at
most $d$, where $d\ge3$. There is a signing
$s:E\to\{-1,+1\}$ such that its signed adjacency matrix satisfies
\[
\|A_s\|<(1+\sqrt2)\sqrt{d-1}.
\]
\end{theorem}

For each $K\subseteq V$, let $\Omega_K=\{-1,+1\}^{E(G[K])}$ be the set of
edge signings of $G[K]$, and define
\[
\mathcal P=\{(K,s):K\subseteq V,\ s\in\Omega_K\}.
\]
In particular, $\mathcal P$ includes $(\varnothing,s_\varnothing)$, where
$s_\varnothing$ is the empty signing of the empty set.
For $(K,s)\in\mathcal P$, let $A_K(s)$ be
the signed adjacency matrix of $G[K]$. We call $(K,s)$ a \emph{good pair} if
$\|A_K(s)\|<(1+\sqrt2)\sqrt{d-1}$.

The proof idea is as follows. We consider the subset
\begin{equation}\label{eq:P-connected}
\mathcal P_{\mathrm{connected}}
=\{(K,s)\in\mathcal P:V\setminus K\ne\varnothing,\ 
G[V\setminus K]\text{ is connected}\}.
\end{equation}
\begin{enumerate}
\item We put a partial order on $\mathcal P$ and take a maximal pair $(K,s)$ from $\mathcal P_{\mathrm{connected}}$ with respect to this order.
\item We prove that this maximal pair misses exactly one vertex from $V$, that is, $K=V\setminus\{u\}$ for some $u\in V$.
\item We prove that the maximal pair $(K,s)$ can be extended to a \emph{good pair} $(V,s')$. Hence, we obtain a good signing on $G$.
\end{enumerate}

\vspace{0.25cm}
\noindent\textcolor{inkpurple}{\textbf{Linear algebra facts used below.}}
We state these facts here and do not prove them.
\begin{enumerate}
\item\label{fact:scalar-schur}
\textcolor{inkpurple}{\textbf{Scalar Schur-complement formula and criterion.}}
If $A$ is invertible and $d$ is a scalar, then
\[
\det\begin{pmatrix}A&b\\ b^{\mathsf T}&d\end{pmatrix}
=\det(A)\bigl(d-b^{\mathsf T}A^{-1}b\bigr).
\]
If $A$ is positive definite, then
\[
\begin{pmatrix}A&b\\ b^{\mathsf T}&d\end{pmatrix}\text{ is positive definite}
\quad\Longleftrightarrow\quad d-b^{\mathsf T}A^{-1}b>0.
\]
\item\label{fact:cofactor}
\textcolor{inkpurple}{\textbf{Cofactor formula for inverse diagonals.}}
If $A$ is invertible and its rows and columns are indexed by $K$, then for
$u\in K$,
\[
(A^{-1})_{uu}=\frac{\det A[K\setminus\{u\}]}{\det A}.
\]
\item\label{fact:inverse-diagonal-lower}
\textcolor{inkpurple}{\textbf{Lower bound for inverse diagonal entries.}}
If $A$ is positive definite, then for every index $u$,
\[
(A^{-1})_{uu}\ge \frac{1}{A_{uu}}.
\]
\end{enumerate}

We are ready to prove Theorem 2.
\begin{proof}[Proof of Theorem 2]
It is enough to prove \hyperref[thm:small-norm-signing]{Theorem~2} for connected graphs $G$ of maximum
degree at most $d$. Indeed, we can apply the result to each connected
component of $G$. The signed adjacency matrix of $G$ is block diagonal,
so its norm is the maximum of the norms of its blocks. Thus, without loss
of generality, we may assume that $G$ is connected.
If $|V|\le1$, the theorem is immediate, so assume $|V|\ge2$. Let $r=(1+\sqrt2)\sqrt{d-1}$.

Define the weight of each pair $(K,s)\in\mathcal P$ by
\[
w_K(s)=
\begin{cases}
\det\!\left(I+\dfrac{A_K(s)}r\right)
\det\!\left(I-\dfrac{A_K(s)}r\right),
&(K,s)\text{ is a \emph{good pair}},\\[1.2ex]
0,&\text{otherwise}.
\end{cases}
\]
For the empty set, the empty signing has weight $w_\varnothing=1$.
Fix $t>1$ and consider the partially ordered set $(\mathcal P,<)$, where
the score of $(K,s)$ is $\Phi(K,s):=t^{|K|}w_K(s)$ and
\[
(K,s)<(K',s')\quad\Longleftrightarrow\quad
\begin{cases}
\Phi(K,s)<\Phi(K',s'), & \text{or}\\[-0.15em]
\Phi(K,s)=\Phi(K',s')\text{ and }|K|<|K'|.
\end{cases}
\]
We take a maximal pair in $\mathcal P_{\mathrm{connected}}$ (as defined in
\eqref{eq:P-connected}) with respect to this order and denote it by $(K,s)$.
The empty pair belongs to $\mathcal P_{\mathrm{connected}}$, since $G$ is
connected, and has weight $1$. To see that it is not maximal, choose a leaf
$u$ of a spanning tree of $G$. Then $G[V\setminus\{u\}]$ is connected, so
$(\{u\},s_{\{u\}})\in\mathcal P_{\mathrm{connected}}$, where
$s_{\{u\}}$ is the empty signing on the singleton. This pair has weight
$1$ and score $t>1$, so it is larger than the empty pair. Hence the maximal
pair $(K,s)$ has positive weight and, by the definition of $w_K(s)$, is a
\emph{good pair}.

Our goal is to prove two claims about this maximal pair:
\begin{enumerate}
\item\label{claim:one-vertex} $K$ misses exactly one vertex from $V$,
that is, $K=V\setminus\{u\}$ for some $u\in V$.
\item\label{claim:extension} It is possible to extend $(K,s)$ to a
\emph{good pair} $(V,s')$. Hence, we obtain a \emph{good signing} on $V$.
\end{enumerate}

We start by proving \hyperref[claim:extension]{claim 2}. Assume that the maximal pair $(K,s)$ misses only
one vertex $u$. We extend $s$ to a signing $s'$ of $G$ by choosing the signs
on the new edges incident to $u$ independently and uniformly from
$\{-1,+1\}$. For each $x\in N_G(u)\cap K$, write
$\varepsilon_x=s'(xu)$ for the random sign on the edge $xu$; these signs are
independent and uniform on $\{-1,+1\}$. Let
\[
Q_K^+=I+\frac{A_K(s)}r,
\qquad
Q_K^-=I-\frac{A_K(s)}r.
\]
Since $(K,s)$ is a \emph{good pair}, the matrices $Q_K^+$ and $Q_K^-$ are
positive definite. Let $c\in\mathbb R^K$ be defined by
$c_x=\varepsilon_x/r$ for $x\in N_G(u)\cap K$ and $c_x=0$ otherwise.
Ordering the vertices of $K$ first and $u$ last, we can write
\[
Q_V^\pm=I\pm\frac{A_V(s')}r
=\begin{pmatrix}Q_K^\pm&\pm c\\\pm c^{\mathsf T}&1\end{pmatrix}.
\]
By the \hyperref[fact:scalar-schur]{scalar Schur-complement fact}, since $Q_K^\pm$ is positive
definite, $Q_V^\pm$ is positive definite exactly when
\[
1-c^{\mathsf T}(Q_K^\pm)^{-1}c>0.
\]
Writing $q_\pm=c^{\mathsf T}(Q_K^\pm)^{-1}c$, both matrices
$Q_V^+$ and $Q_V^-$ are positive definite exactly when
both $q_\pm<1$. We will show that this happens with positive probability,
yielding an extension $s'$ with $\|A_V(s')\|<r$. Thus $(V,s')$ is a
\emph{good pair} and $s'$ a \emph{good signing} of $G$.

We calculate the expectation of $q_++q_-$ over the random signs
$\varepsilon_x$. Since $c_x=\varepsilon_x/r$ on $N_G(u)\cap K$ and is zero
otherwise, independence and $\mathbb E[\varepsilon_x]=0$ give
\[
\mathbb E[c_xc_y]=
\begin{cases}
r^{-2},&x=y\in N_G(u)\cap K,\\
0,&\text{otherwise}.
\end{cases}
\]
Therefore,
\begin{equation}\label{eq:extension-expectation}
\small
\mathbb E[q_++q_-]
=\mathbb E\!\left[c^{\mathsf T}\left((Q_K^+)^{-1}+(Q_K^-)^{-1}\right)c\right]
=\frac1{r^2}\sum_{x\in N_G(u)\cap K}
\left(\bigl((Q_K^+)^{-1}\bigr)_{xx}+\bigl((Q_K^-)^{-1}\bigr)_{xx}\right).
\end{equation}
By the \hyperref[fact:cofactor]{cofactor formula for inverse diagonals} stated above, for each
$x\in N_G(u)\cap K$ and each choice of sign,
\begin{equation}\label{eq:inverse-diagonal-cofactor}
\bigl((Q_K^\pm)^{-1}\bigr)_{xx}
=\frac{\det Q_K^\pm[K\setminus\{x\}]}
{\det Q_K^\pm}.
\end{equation}
Since $x$ is adjacent to the remaining vertex $u$, the complement of
$K\setminus\{x\}$ is $\{u,x\}$, which is connected. Thus the restricted
signing gives another pair in $\mathcal P_{\mathrm{connected}}$. It is also
a good pair: the matrices $Q_{K\setminus\{x\}}^\pm$ are principal
submatrices of the positive definite matrices $Q_K^\pm$, and hence are
positive definite.
By maximality,
\begin{equation}\label{eq:maximality-deletion}
t^{|K|}w_K(s)
\ge t^{|K|-1}
w_{K\setminus\{x\}}\!\left(s\big|_{E(G[K\setminus\{x\}])}\right).
\end{equation}
For lighter notation, write
$w_{K\setminus\{x\}}:=w_{K\setminus\{x\}}\!\left(s\big|_{E(G[K\setminus\{x\}])}\right)$.
By the definition of the weights and equation~\eqref{eq:inverse-diagonal-cofactor},
\[
\small
\frac{w_{K\setminus\{x\}}}{w_K(s)}
=\frac{\det Q_K^+[K\setminus\{x\}]\,\det Q_K^-[K\setminus\{x\}]}
{\det Q_K^+\,\det Q_K^-}
=\bigl((Q_K^+)^{-1}\bigr)_{xx}\bigl((Q_K^-)^{-1}\bigr)_{xx}.
\]
Combining this with the maximality inequality~\eqref{eq:maximality-deletion},
we conclude that
\[
\bigl((Q_K^+)^{-1}\bigr)_{xx}
\bigl((Q_K^-)^{-1}\bigr)_{xx}\le t.
\]
By the \hyperref[fact:inverse-diagonal-lower]{lower bound for inverse diagonal entries} stated above, and since the
diagonal entries of $Q_K^+$ and $Q_K^-$ are all $1$, set
\[
a_x:=\bigl((Q_K^+)^{-1}\bigr)_{xx}\ge1,
\qquad
b_x:=\bigl((Q_K^-)^{-1}\bigr)_{xx}\ge1.
\]
For any $a,b\ge1$,
\[
\begin{aligned}
0&\le(a-1)(b-1)=ab-a-b+1,
\end{aligned}
\]
which implies $a+b\le1+ab$. Substituting $a=a_x$ and $b=b_x$, and using
$a_xb_x\le t$, we conclude that
\[
\bigl((Q_K^+)^{-1}\bigr)_{xx}
+\bigl((Q_K^-)^{-1}\bigr)_{xx}\le1+t.
\]
Returning to~\eqref{eq:extension-expectation}, we obtain
\[
\mathbb E[q_++q_-]
\le\frac{(1+t)|N_G(u)\cap K|}{r^2}
\le\frac{(1+t)d}{r^2}.
\]
Choose $t=1+\sqrt2$, so that $r=t\sqrt{d-1}$. Then
\[
\mathbb E[q_++q_-]
\le\frac{(1+t)d}{t^2(d-1)}
=\left(1-\frac1t\right)\frac{d}{d-1}
\le\frac32\left(1-\frac1t\right)<1.
\]
Therefore some choice of the new signs satisfies $q_++q_-<1$. Since
$q_+,q_-\ge0$, this gives $q_\pm<1$ for both choices of sign, and hence
both $Q_V^+$ and $Q_V^-$ are positive definite.

All that remains to prove is \hyperref[claim:one-vertex]{claim 1},
namely that $K=V\setminus\{u\}$ for some $u\in V$.
We prove this by contradiction. Assume that $R=V\setminus K$ satisfies
$|R|>1$. Since $(K,s)\in\mathcal P_{\mathrm{connected}}$, the graph $G[R]$
is connected. Choose a spanning tree of $G[R]$ and let $u$ be a leaf. It has
a neighbor still in $R$, so it has at most
$d-1$ neighbors in $K$. Also, $R\setminus\{u\}$ is nonempty and connected,
so for any signing $s'$ extending $s$ to $K\cup\{u\}$, the pair
$(K\cup\{u\},s')$ belongs to $\mathcal P_{\mathrm{connected}}$.
Similarly to the argument for claim 2, we extend $(K,s)$ to a pair
$(K\cup\{u\},s')$ by choosing the signs on the new edges incident to $u$
independently and uniformly from $\{-1,+1\}$. We will show that some choice
of these signs gives a pair whose score is at least that of $(K,s)$,
contradicting maximality.

As before, we denote by $c\in\mathbb R^K$ the vector with
$c_x=s'(xu)/r$ for $x\in N_G(u)\cap K$ and $c_x=0$ otherwise, and set
$q_\pm=c^{\mathsf T}(Q_K^\pm)^{-1}c$. As before, the extension is a
\emph{good pair} exactly when $q_\pm<1$ for both signs.
For each $x\in N_G(u)\cap K$, the complement of $K\setminus\{x\}$ is
$R\cup\{x\}$, which is connected because $x$ is adjacent to $u\in R$.
Thus the same maximality argument as above gives
\[
\bigl((Q_K^+)^{-1}\bigr)_{xx}
+\bigl((Q_K^-)^{-1}\bigr)_{xx}\le1+t.
\]
Returning to the expectation formula~\eqref{eq:extension-expectation} and
using this time $|N_G(u)\cap K|\le d-1$, we obtain
\[
\mathbb E[q_++q_-]
\le\frac{(1+t)(d-1)}{r^2}
=\frac{1+t}{t^2}
=1-\frac1t,
\]
Hence some choice of the new signs satisfies $q_++q_-\le1-1/t<1$.
For this choice, the extension is a \emph{good pair}.

For this good extension, the \hyperref[fact:scalar-schur]{scalar
Schur-complement determinant formula} gives
\[
\begin{aligned}
w_{K\cup\{u\}}(s')
&=\det Q_{K\cup\{u\}}^+\det Q_{K\cup\{u\}}^-\\
&=\det Q_K^+\det Q_K^-(1-q_+)(1-q_-)\\
&=w_K(s)(1-q_+)(1-q_-)\\
&\ge\frac1t w_K(s).
\end{aligned}
\]
Indeed, $(1-q_+)(1-q_-)=1-(q_++q_-)+q_+q_-\ge1-(q_++q_-)\ge1/t$.
Consequently,
\[
t^{|K\cup\{u\}|}w_{K\cup\{u\}}(s')
\ge t^{|K|}w_K(s).
\]
Since $|K\cup\{u\}|>|K|$, the new pair is larger in the order, a
contradiction to maximality.
\end{proof}

Using the same technique, we also recover the celebrated bipartite signing
bound of Marcus, Spielman, and Srivastava~\cite{MSS}.
In the bipartite case, symmetry gives $q_+=q_-$. Thus a bound on $q_+$ controls both Schur-complement conditions at once.
\begin{corollary}[Bipartite signing bound]\label{thm:bipartite-signing}
For every finite bipartite graph $G=(V,E)$ of maximum degree at most
$d\ge3$, there is an edge signing $s$ such that
\[
\|A_s\|<2\sqrt{d-1}.
\]
\end{corollary}
\begin{proof}[Proof of Corollary~3]
As in the proof of Theorem~2, we may assume that $G$ is connected and
has at least two vertices. Set
\[
r=2\sqrt{d-1},\qquad t=4.
\]
Throughout this proof, call $(K,s)$ a \emph{good pair} if
$\|A_K(s)\|<r$, and define its weight by
\[
w_K(s)=
\begin{cases}
\det\!\left(I+\dfrac{A_K(s)}r\right)
\det\!\left(I-\dfrac{A_K(s)}r\right), &(K,s)\text{ is a good pair},\\[1.2ex]
0, &\text{otherwise}.
\end{cases}
\]
Use the same partial-order construction as in Theorem~2, with score
$\Phi(K,s)=t^{|K|}w_K(s)$ using this weight. Let $(K,s)$ be a maximal pair
in $\mathcal P_{\mathrm{connected}}$, and write $R=V\setminus K$.
We prove the same two claims:
\begin{enumerate}
\item $K$ misses exactly one vertex, that is, $R=\{u\}$ for some $u\in V$.
\item If $R=\{u\}$, then $(K,s)$ extends to a \emph{good pair}
$(V,s')$.
\end{enumerate}

Fix $u\in R$. As in the proof of Theorem~2, choose the signs on the
edges from $u$ to $K$ independently and uniformly. Let
$c_x=s'(xu)/r$ for $x\in N_G(u)\cap K$, and let $c_x=0$ otherwise.
Recall that
\[
Q_K^\pm=I\pm\frac{A_K(s)}r,
\qquad
q_\pm=c^{\mathsf T}(Q_K^\pm)^{-1}c.
\]
By the same Schur-complement argument, the extension is a \emph{good pair}
exactly when both $q_+<1$ and $q_-<1$.

Fix a bipartition $V=V_+\sqcup V_-$ of $G$. For each $K\subseteq V$,
let $J_K$ be the diagonal matrix with entry $+1$ on $K\cap V_+$ and
$-1$ on $K\cap V_-$. Since every edge joins the two parts,
$J_KA_K(s)J_K=-A_K(s)$. Hence
\[
Q_K^-=J_KQ_K^+J_K,
(Q_K^-)^{-1}=J_K(Q_K^+)^{-1}J_K.
\]
All neighbors of $u$ in $K$ lie in the same part of the fixed
bipartition, so $J_Kc=\pm c$. It follows that
\[
q_-
=c^{\mathsf T}(Q_K^-)^{-1}c
=c^{\mathsf T}J_K(Q_K^+)^{-1}J_Kc
=c^{\mathsf T}(Q_K^+)^{-1}c
=q_+.
\]
Write $q$ for their common value. Fix $x\in N_G(u)\cap K$. Write $a_x$ and $b_x$ for the two inverse
diagonal entries, as in the proof of Theorem~2:
\[
a_x=\bigl((Q_K^+)^{-1}\bigr)_{xx},
\qquad
b_x=\bigl((Q_K^-)^{-1}\bigr)_{xx}.
\]
Just as in that proof, the cofactor formula and maximality give
\[
a_xb_x=\frac{w_{K\setminus\{x\}}(s|_{K\setminus\{x\}})}{w_K(s)}\le t.
\]
Since $(Q_K^-)^{-1}=J(Q_K^+)^{-1}J$ and conjugation by $J$ leaves
each diagonal entry unchanged, we have $a_x=b_x$. Thus
\[
a_x^2=a_xb_x\le t,
\qquad\text{so}\qquad
a_x=b_x\le\sqrt t=2.
\]
Averaging over the new signs exactly as in the proof of Theorem~2, we
obtain
\[
\mathbb E[q]
=\frac{1}{r^2}\sum_{x\in N_G(u)\cap K}a_x
\le\frac{\sqrt t\,|N_G(u)\cap K|}{r^2}
=\frac{|N_G(u)\cap K|}{2(d-1)}.
\]

We first prove claim~2. If $R=\{u\}$, then $u$ has at most $d$
neighbors in $K$. Hence
\[
\mathbb E[q]\le\frac{d}{2(d-1)}<1,
\]
since $d\ge3$. Some choice of the new signs therefore gives $q<1$.
Since $q_+=q_-=q$, both Schur-complement conditions hold, and this
choice extends $(K,s)$ to a \emph{good pair} $(V,s')$.

It remains to prove claim~1. Suppose $|R|>1$. As in Theorem~2, choose
a leaf $u$ of a spanning tree of $G[R]$. Then $u$ has at most $d-1$
neighbors in $K$, and $G[R\setminus\{u\}]$ is nonempty and connected.
This time,
\[
\mathbb E[q]\le\frac{d-1}{2(d-1)}
=\frac12
=1-\frac1{\sqrt t}.
\]
Choose signs for which $q\le 1-1/\sqrt t$. The extension is a
\emph{good pair}, and the Schur-complement determinant formula gives
\[
w_{K\cup\{u\}}(s')
=w_K(s)(1-q_+)(1-q_-)
=w_K(s)(1-q)^2
\ge \frac1t\,w_K(s).
\]
Thus the extended pair has at least the same score and a larger vertex
set, contradicting maximality. Hence $R$ consists of one vertex, and
claim~2 gives a signing $s'$ of $G$ with
$\|A_{s'}\|<2\sqrt{d-1}$.
\end{proof}

\phantomsection
\addcontentsline{toc}{section}{References}
\bibliographystyle{plain}
\bibliography{bilu-linial-signing-counterexamples-and-bounds}

@article{BiluLinial,
  author  = {Bilu, Yonatan and Linial, Nathan},
  title   = {Lifts, discrepancy and nearly optimal spectral gap},
  journal = {Combinatorica},
  volume  = {26},
  number  = {5},
  pages   = {495--519},
  year    = {2006},
  note    = {\href{https://doi.org/10.1007/s00493-006-0029-7}{doi:10.1007/s00493-006-0029-7}}
}

@article{MSS,
  author  = {Marcus, Adam W. and Spielman, Daniel A. and Srivastava, Nikhil},
  title   = {Interlacing Families {I}: Bipartite Ramanujan Graphs of All Degrees},
  journal = {Annals of Mathematics},
  volume  = {182},
  number  = {1},
  pages   = {307--325},
  year    = {2015},
  note    = {\href{https://doi.org/10.4007/annals.2015.182.1.7}{doi:10.4007/annals.2015.182.1.7}; \href{https://arxiv.org/abs/1304.4132}{arXiv:1304.4132}}
}

@article{RS2017,
  author  = {Ravichandran, Mohan and Srivastava, Nikhil},
  title   = {Asymptotically Optimal Multi-Paving},
  journal = {International Mathematics Research Notices},
  volume  = {2021},
  number  = {14},
  pages   = {10908--10940},
  year    = {2021},
  note    = {\href{https://doi.org/10.1093/imrn/rnz111}{doi:10.1093/imrn/rnz111}; \href{https://arxiv.org/abs/1706.03737}{arXiv:1706.03737}}
}

@misc{Huang2026,
  author       = {Huang, Zejun},
  title        = {A $\sqrt{2}$-Approximation to the Bilu--Linial Conjecture},
  year         = {2026},
  howpublished = {arXiv:2609.17551},
  note         = {\href{https://arxiv.org/abs/2609.17551}{arXiv:2609.17551}}
}

@misc{LinZhou2026,
  author       = {Lin, Fangfang and Zhou, Hong},
  title        = {Improved Bounds for the Bilu--Linial Conjecture via Spectral Recovery from Mixed Determinantal Polynomials},
  year         = {2026},
  howpublished = {arXiv:2609.15715},
  note         = {\href{https://arxiv.org/abs/2609.15715}{arXiv:2609.15715}}
}

@misc{JSS2026,
  author       = {Jadbabaie, Ali and Saberi, Amin and Sra, Suvrit},
  title        = {Kadison--Singer Partitions and Bilu--Linial Graph Signings in Polynomial Time},
  year         = {2026},
  howpublished = {arXiv:2609.23855},
  note         = {\href{https://arxiv.org/abs/2609.23855}{arXiv:2609.23855}}
}

@misc{Xu2026,
  author       = {Xu, Zhiqiang},
  title        = {A $3$-Regular Counterexample to the Bilu--Linial Signing Conjecture},
  year         = {2026},
  howpublished = {arXiv:2609.15591},
  note         = {\href{https://arxiv.org/abs/2609.15591}{arXiv:2609.15591}}
}

\end{document}